\documentclass[conference]{IEEEtran}
\IEEEoverridecommandlockouts
\usepackage[space]{cite}
\usepackage{amsmath,amssymb,amsfonts}
\usepackage{algorithmic}
\usepackage{graphicx}
\usepackage{textcomp}
\usepackage[dvipsnames]{xcolor}
\usepackage{caption}
\usepackage{subcaption}
\usepackage{hyperref}
\usepackage{tikz}
\usepackage{pgfplots}
\usepackage[squaren]{SIunits}
\usepackage{pgfplots}
\pgfplotsset{compat=1.15}
\usepackage{mathrsfs}
\usetikzlibrary{arrows}
\newcommand{\timeidx}{t}
\newcommand{\defeq}{:=}
\newcommand{\starttime}{t_{\rm s}}
\newcommand{\durationload}{T_{\rm a}}
\newcommand{\windpower}[1]{P^{(\rm w)}_{#1}}
\newcommand{\appliance}[1]{P^{(\rm a)}_{#1}} 
\newcommand{\dishwasher}[1]{P^{(\rm dw)}_{#1}} 
\newcommand{\household}[1]{P^{(\rm house)}_{#1}} 
\newcommand{\batterylevel}[1]{E^{(\rm b)}_{#1}} 
\newcommand{\windspeed}{V} 
\newcommand{\battcap}{E_{bat}}
\newcommand{\fractiongeneric}{\rho}
\newcommand{\fractionfmi}[1]{\rho^{(#1)}}
\newcommand{\stationindex}{i} 
\newcommand{\houridx}{i} 

\pgfplotsset{compat=newest}
\hypersetup{
    colorlinks=true,
    linkcolor=blue,
    filecolor=magenta,      
    urlcolor=cyan,
    pdftitle={Overleaf Example},
    pdfpagemode=FullScreen,
    }

\def\BibTeX{{\rm B\kern-.05em{\sc i\kern-.025em b}\kern-.08em
    T\kern-.1667em\lower.7ex\hbox{E}\kern-.125emX}}

\begin{document}

\title{Wind to start the washing machine? \\ High-Resolution Wind Atlas for Finland}

\author{\IEEEauthorblockN{1\textsuperscript{st} Xu Yang}
\IEEEauthorblockA{\textit{Dept. of Computer Science} \\
\textit{Aalto University}\\
Espoo, Finland \\
xu.1.yang@aalto.fi}
\and
\IEEEauthorblockN{2\textsuperscript{nd} Yu Tian}
\IEEEauthorblockA{\textit{Dept. of Computer Science} \\
\textit{Aalto University}\\
Helsinki, Finland \\
yu.tian@aalto.fi}
\and
\IEEEauthorblockN{3\textsuperscript{rd} Irene Schicker}
\IEEEauthorblockA{\textit{Dept. of Postprocessing} \\
	\textit{GeoSphere Austria}\\
	Vienna, Austria \\
	irene.schicker@geosphere.at}
\and
\IEEEauthorblockN{4\textsuperscript{th} Alexander Jung}
\IEEEauthorblockA{\textit{Dept. of Computer Science} \\
\textit{Aalto University}\\
Helsinki, Finland \\
alex.jung@aalto.fi}
}

\maketitle

\begin{abstract}
The current fossil fuel and climate crisis has led to an increased demand for renewable energy sources, such as wind power. In northern Europe, the efficient use of wind power is crucial for achieving carbon neutrality. To assess the potential of wind energy for private households in Finland, we have conducted a high spatiotemporal resolution analysis. Our 
main contribution is a wind power map of Finland that indicates the availability 
of wind power for given load profiles. As a representative example of power load, 
we consider the load profile of a household appliance. We compare this load profile against 
the wind power available nearby the weather stations of the Finnish meteorological institute. 
\end{abstract}


\section{Introduction}

The efficient use of renewable energy is a key component for the green transition of current 
fossil-fuel based industries to achieve carbon neutrality \cite{carbonneutrality}. Wind energy is a 
main source for renewable energy in Nordic countries such as Finland, especially during winter 
months \cite{energy2022}. According to the up-to-date statistics(2022), wind power has covered 14,1\% of the Finnish
electricity consumption\cite{statistics2022}.




A key challenge in the efficiently using of wind energy is its spatio-temporal variability. 
In general, wind energy is available at locations and during times, which may not always match the end-user's needs (household appliance). The transfer of wind energy across space 
and time requires well-designed power grids and efficient storage facilities such as batteries \cite{qi2015joint, willis2018wind}. 

Since the transmission of wind energy incurs losses \cite{morales2012transmission}, it is beneficial 
to consume the wind energy near its production sites.This highlights the importance of site selection for wind power plants, which often requires corresponding estimation of wind production. Previous studies \cite{wais2017two,baseer2017wind,celik2004statistical,azad2014statistical} have found that 
Weibull distribution is a useful tool for the evaluation of wind resources for chosen sites, due to its ability to 
properly fit wind data as a probability distribution function.

Historical wind speed observation data of Malaysia has been used to evaluate wind power density and inform wind power plant site selection in \cite{masseran2012analysis}.
The authors of \cite{bennui2007site} combined geographic information systems with multi-criteria decision making 
to optimize the site selection in Thailand. These works use a long-term average perspective, based on the 
statistics of wind speed and power production. In contrast, we are interested in the short-term availability of 
wind power to operate household appliances, equipped with modest battery capacity. 

Our work is most closely related to recent efforts in generating various instances of a 
wind power atlas \cite{newa, era5}. In particular, the NEWA and ERA5 reanalysis data sets have been proposed 
for wind analysis fields with a temporal resolution of $30$ and $60$ minute intervals, respectively, and spatial resolutions of $3$ and $30$ km, respectively. 
In contrast to these existing works, our approach uses a higher temporal resolution with $10$ minute 
intervals. On the other hand, our approach focuses on local wind power generation dictated by the locations of 
weather stations operated by the Finnish Meteorological Institute (FMI). 



{\bf Contribution.} This paper provides the results of an exploratory data analysis using freely available 
weather data provided by FMI. The aim of this analysis is to generate a high-resolution wind-power map for Finland. For each FMI station we determine the fraction of the year 2021 during which an 
appliance with a given power profile could be powered solely from wind power. 

{\bf Notation.} We denote the first $n$ natural numbers starting with $0$ as $[n] \defeq \{0,\ldots,n-1\}$. $V$ for wind speed, $P$ for wind power, and $E$ for energy are used throughout the paper. Other symbols are defined as required.

\section{Problem Setting}

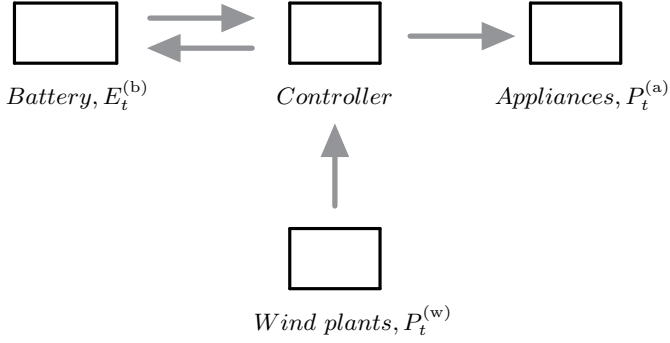
\begin{figure} 
	
	\definecolor{zzttqq}{rgb}{0.6,0.2,0}
	\begin{tikzpicture}[line cap=round,line join=round,>=triangle 45,x=2cm,y=2cm]
		\clip(-0.65,-0.2) rectangle (6,2.3);
        \draw [line width=1.2pt] (-0.58,1.7)-- (0.1,1.7);
        \draw [line width=1.2pt] (0.1,2.1)-- (0.1,1.7);
        \draw [line width=1.2pt] (0.1,2.1)-- (-0.58,2.1);
        \draw [line width=1.2pt] (-0.58,2.1)-- (-0.58,1.7);
        \draw (-0.7,1.65) node[anchor=north west, font=\fontsize{9}{0}\selectfont] {$Battery, \batterylevel{\timeidx}$};

		\draw [->,line width=2pt,color=Gray] (1,1.8) -- (0.3,1.8);

		\draw [line width=1.2pt] (1.25,1.7)-- (1.85,1.7);
		\draw [line width=1.2pt] (1.85,1.7)-- (1.85,2.1);
		\draw [line width=1.2pt] (1.85,2.1)-- (1.25,2.1);
		\draw [line width=1.2pt] (1.25,2.1)-- (1.25,1.7);
		
		\draw [->,line width=2pt,color=Gray] (2.05,1.88) -- (2.75,1.88);
		\draw [->,line width=2pt,color=Gray] (1.55,0.75) -- (1.55,1.3);

	    \draw [line width=1.2pt] (2.85,2.1)-- (2.85,1.7);
	    \draw [line width=1.2pt] (2.85,1.7)-- (3.45,1.7);
	    \draw [line width=1.2pt] (3.45,1.7)-- (3.45,2.1);
	    \draw [line width=1.2pt] (2.85,2.1)-- (3.45,2.1);

	    \draw [line width=1.2pt] (1.25,0.6)-- (1.25,0.2);
	    \draw [line width=1.2pt] (1.25,0.6)-- (1.85,0.6);
	    \draw [line width=1.2pt] (1.85,0.6)-- (1.85,0.2);
	    \draw [line width=1.2pt] (1.85,0.2)-- (1.25,0.2);

		\draw (1.1,1.6) node[anchor=north west, font=\fontsize{9}{0}\selectfont] {$Controller$};
		\draw (0.95,0.15) node[anchor=north west, font=\fontsize{9}{0}\selectfont] {$Wind$ $plants, \windpower{\timeidx}$};
		\draw (2.55,1.65) node[anchor=north west, font=\fontsize{9}{0}\selectfont] {$Appliances, \appliance{\timeidx}$};
		\draw [->,line width=2pt,color=Gray] (0.32,2) -- (1.02,2);
		
	\end{tikzpicture}
	\caption{A wind power plant generates the power $\windpower{\timeidx}$ at discrete-time $\timeidx$ which 
		is used to serve a load with prescribed power profile $\appliance{\timeidx}$. The surplus (if any) power is 
		used to load a battery whose energy level is $\batterylevel{\timeidx}$. 
		\label{fig_simple_setup}
	}
\end{figure}

We consider the simple power system depicted in Figure \ref{fig_simple_setup} during 
discrete time instants $\timeidx=0,1,\ldots$. The absolute time difference between any two 
consecutive time instants $\timeidx$ and $\timeidx\!+\!1$ is $\Delta \timeidx =10 \mbox{min}$. 
The system includes a wind power plant that delivers the power $\windpower{\timeidx}$ at 
time instant $\timeidx$. We consider a wind power plant of type  \texttt{Nordex N100/25000} \url{https://www.thewindpower.net/turbine_en_224_nordex_n100-2500.php}  that is mounted at a height of $100$ m. 

The system in Figure \ref{fig_simple_setup} also includes a load that is characterized by a 
power profile $\appliance{\timeidx'}$ for time instants $ \timeidx' \in [\durationload]$. 
The total duration (in absolute time) of the load profile is $\durationload \cdot 1 {\rm min}$. 
An example for the load is a household appliance such as a dishwasher or washing 
machine (see Figure \ref{fig:dishwasher}). The power profile of the load has finite 
support of $\durationload$ time instants, $$\appliance{\timeidx'}=0 \mbox{ for }\timeidx' \notin \{0,1,\ldots,\durationload \}.$$

Consider some candidate time instant $\starttime \in [365 \cdot 24 \cdot 6]$ 
during the year $2021$. Starting at $\starttime$ we try to run the load $\appliance{\timeidx}$. 
We assume the battery is empty when starting the load and ignore any power leakage, 
\begin{align}
\batterylevel{\starttime} &= 0,  \nonumber \\ 
\batterylevel{\timeidx\!+\!1} &= \min\big\{ \batterylevel{\timeidx}\!+\! \big( \windpower{\timeidx}\!-\!\appliance{\timeidx\!-\!\starttime}\big) \Delta \timeidx, \battcap \big\} \mbox{ for } \timeidx\!>\!\starttime. 
\label{equ:work}
\end{align} 
We define the candidate starting time $\starttime$ as {\bf suitable} if $\batterylevel{\timeidx}  \geq 0$ 
for $\timeidx \in \{\starttime,\starttime+1,\ldots,\starttime+ \durationload\}$. The useful annual
fraction of the year $2021$ is defined as 
\begin{equation}
	\fractiongeneric \defeq  \frac{\big| \big\{ \starttime \in [365 \cdot 24 \cdot 6]: \starttime \mbox{ is {\bf suitable}} \big\}  \big|}{365 \cdot 24 \cdot 6}
	\label{equ:frac}
\end{equation} 
Note that the fraction $\fractiongeneric$ depends on the battery capacity $\battcap$, the load profile $\appliance{\timeidx'}$ and the 
available wind power $\windpower{\timeidx'}$. Section \ref{sec_wind_power_map} discusses different choices 
for the batter capacity $\battcap$ and and load profiles $\appliance{\timeidx'}$ that are used to determine 
the useful annual fractions $\fractionfmi{\stationindex}$ nearby FMI weather stations, indexed by $\stationindex \in \{1,2,\ldots\}$.  


\section{Method}
\label{sec_wind_power_map} 

To determine the useful factions $\fractionfmi{\stationindex}$ nearby FMI station $\stationindex=1,2,\ldots$, 
we estimate the available wind power from the wind speed observations at a height of $10$ m. 
Section \ref{sec_data_preprocessing} explains the pre-processing of raw weather observations, including 
the imputation of missing observations, extropalation of wind speed at $10$ m to $100$ m, and  the estimation of wind 
power generation at $100$ m. Section \ref{sec_power_load} discusses representative power load profiles that we will 
use for computing the useful fractions and generating the wind power atlas.

\subsection{Pre-Processing Wind Observations}
\label{sec_data_preprocessing} 

We downloaded wind speed observations $\windspeed_{\timeidx}$ (height $= 10$ m) at time instants $\timeidx= \{0,1,\ldots,365 \cdot 24 \cdot 6 \}$ 
during year $2021$ at different FMI stations from the web interface \url{https://en.ilmatieteenlaitos.fi/open-data}. 
In a next step we excluded any weather station for which more than $3\%$ wind observations were missing. 
The resulting $165$ FMI weather stations, indexed $\stationindex=1,2,\ldots,165$, are then used to construct 
the wind atlas in Section \ref{sec:atlas}. 

{\bf Data Imputation.}  For some weather stations, the wind speed observations are missing (we also treat negative wind speed values as missing values) at some time instants. 
We choose to impute missing wind observations via linear temporal interpolation\cite{federer1996intercomparison}. 
If we denote $\timeidx_{m}$ and $\timeidx_{n}$ the time instants just before and after the time instants 
of missing observations, 
\begin{equation}
     \widehat{\windspeed}_{\timeidx}\!=\!\frac{\windspeed_{\timeidx_m}\!-\!\windspeed_{\timeidx_n}}{\timeidx_m\!-\!\timeidx_n }  \cdot\left(\timeidx\!-\!\timeidx_m\right)\!+\!\windspeed_{\timeidx_m} \mbox{ for } \timeidx\!\in\!\{ \timeidx_{m}\!+\!1,\ldots,\timeidx_{n}\!-\!1 \}. 
\end{equation}


{\bf From Wind Speed to Wind Power.}: We use the ``1/7 wind power law'' \cite{de2005venerable} to  extrapolate the 
wind speed measured at the height of $10$ m to the expected wind speed at $100$ m (the wind turbine hub height):
\begin{equation}
\widehat{\windspeed}_{\text {100}}=\windspeed_{\text {10 }} \cdot\left(\frac{100}{10}\right)^\alpha
\end{equation}
with $$\alpha=1/7.$$
The estimated wind speed $\widehat{\windspeed}_{\text {100}}$ at each time instant $\timeidx$ 
is then combined with the power curve of the turbine  \texttt{Nordex N100/25000} to 
obtain an estimate for the generated power $\windpower{\timeidx}$. 
In particular,  the estimated power $\widehat{P}$ delivered by the turbine for an estimated wind 
speed in the range $\windspeed_j <\widehat{\windspeed}_{100} < \windspeed_{j'}$ 
is 
\begin{equation}
    \widehat{P}=\frac{P_{j}-P_{j'}}{\windspeed_{j}-\windspeed_{j'}} \cdot\left(\widehat{\windspeed}_{100}-\windspeed_{j}\right)+P_{j} 
\end{equation}
Here, $P_{j}$ and $P_{j'}$ denote, respectively, the nominal wind power delivered at wind speeds $\windspeed_j$ and $\windspeed_{j'}$.

\subsection{Power Load Profiles}
\label{sec_power_load} 

We use the open dataset \url{https://www.kaggle.com/datasets/uciml/electric-power-consumption-data-set} to 
construct two prototype load profiles $\appliance{\timeidx'}$. The first load profile $\dishwasher{\timeidx}$ 
corresponds to a single dishwasher. The second load profile $\household{\timeidx'}$ 
represents an entire single-family household. 
	
\begin{figure}
\begin{tikzpicture}
\begin{axis}[
	ylabel=$\dishwasher{\timeidx'} \mbox{ [Watt]}$, 
	xlabel=$\timeidx' \mbox{ [min]}$,
	ytick={0,500,1000,1500,2000,2500},
	xtick={0,10,20,30,40,50,60,70,80},
	grid=both,
	minor grid style={gray!50},
	major grid style={gray!50},
	width=0.75\linewidth,
	no marks]
\addplot[line width=1pt,solid,color=NavyBlue] %
	table[x=time,y=power,col sep=comma]{dishwasher.csv};
\end{axis},
\end{tikzpicture}
\caption{The load profile of a dishwasher with time intervals of  $1$ minute.}
\label{fig:dishwasher}
\end{figure}
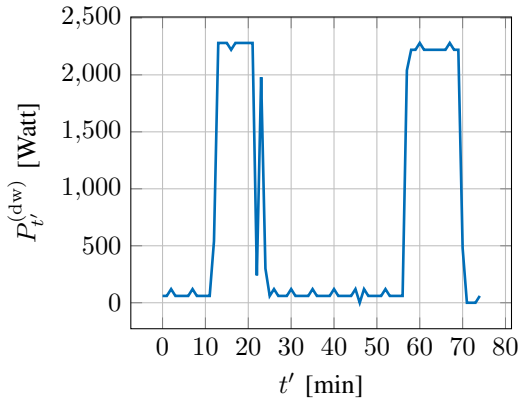

{\bf Single Appliance.}
\label{sec:dish}
We extracted the load power profile $\dishwasher{\timeidx}$  (Figure \ref{fig:dishwasher} ) of a dishwasher from the dataset \cite{paganelli2015appliance}. 

The overall duration of the dishwasher process is $75$ minutes with time intervals of  $1$ minute. To align the different time intervals of dishwasher load profile and wind speed records from FMI, wind speed is assumed to be constant during each 10 min time interval. Based on this assumption, we implemented  (\ref{equ:work})  on a 1-min interval to run the corresponding numerical experiments.

\begin{figure}
	\begin{tikzpicture}
		\begin{axis}[
			ylabel=$\household{\timeidx'} \mbox{ [Watt]}$,
			xlabel=$\timeidx' / 6$,
			ytick={0,500,1000,1500,2000,2500,3000,3500},
			xtick = {0,23,47,71,95,119,143},
			xticklabels={0,4,8,12,16,20,24},
			grid=both,
			minor grid style={gray!50},
			major grid style={gray!50},
			width=0.75\linewidth,
			no marks]
			\addplot[line width=1pt,solid,color=NavyBlue] %
			table[x=time,y=power,col sep=comma]{house.csv};
		\end{axis},
	\end{tikzpicture}
	\caption{Load profile of a representative household during an entire day which corresponds to 
		a duration of $\durationload = 24 \cdot 6$ time instants ($10$ minute intervals).}
	\label{fig:house}
\end{figure}
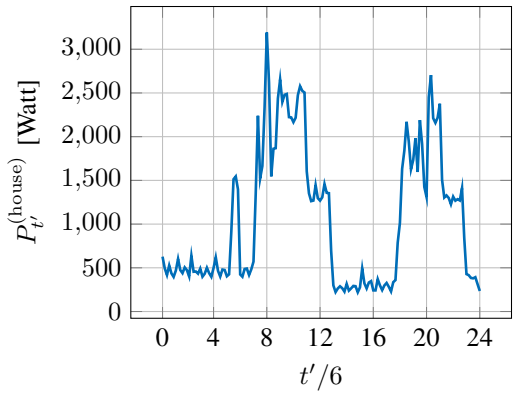

{\bf Entire Household.}
Figure \ref{fig:house} depicts a representative load profile $\household{\timeidx'}$ for a private household 
during an entire day \cite{paganelli2015appliance}. The duration of the profile $\household{\timeidx'}$ is 
$\durationload = 24 \cdot 6$ time instants. 


\section{Wind Power Atlas for Finland} 
\label{sec:atlas}

Numerical experiments were then conducted to explore fine-resolution temporal and spatial variations of the wind power resource in Finland via virtually running 
this dishwasher process according to Eq. (\ref{equ:work}) and Eq. (\ref{equ:frac}). Table \ref{tab:dishwasher} shows useful annual fractions at FMI weather stations, given different battery capacities. The useful fractions increase as the battery capacity increases until a certain value (between $800$ and $1000 {\rm Wh}$) 
is reached. Increasing battery capacity beyond this value has no effect on the resulting useful annual fractions. 

\begin{figure}[htb]
   \center{\includegraphics[width=60mm]
	       {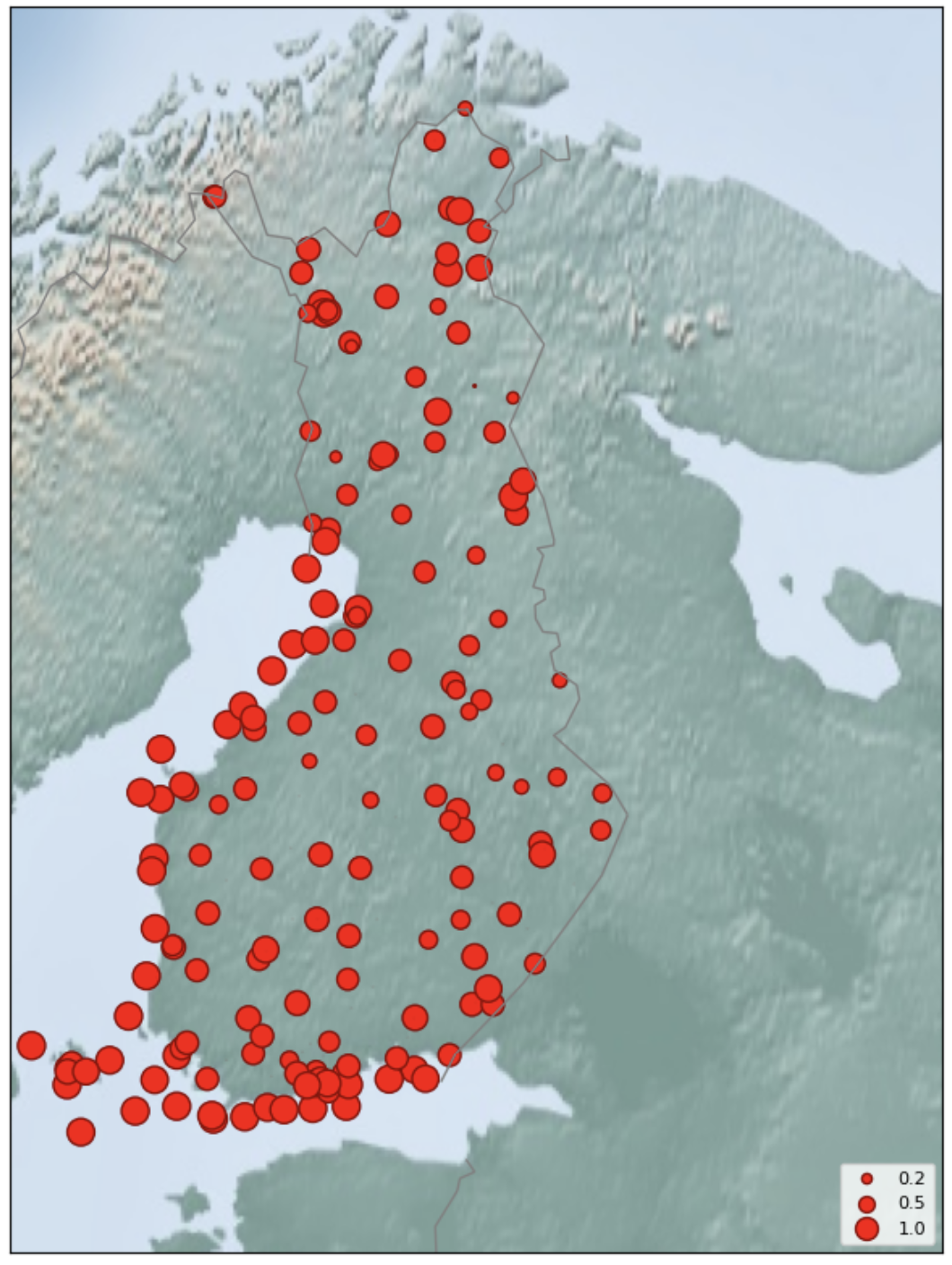}}
	 \caption{Spatial distribution of useful annual fractions $\fractionfmi{\stationindex}$ of $2021$ during which the 
	 	load profile $\dishwasher{\timeidx'}$ (see Figure \ref{fig:dishwasher}) could have been powered only from 
	 	wind power (using a battery capacity $\battcap=1000 {\rm Wh}$). Each dot represents an FMI weather station, 
	 	indexed by $\stationindex=1,\ldots,165$ and its radius is scaled by $\fractionfmi{\stationindex}$. }
	 \label{fig:heat}
\end{figure}

Figure \ref{fig:heat} is a map of Finland, in the location of FMI weather station $\stationindex$, a red dot is 
added as a marker; the marker size represents the corresponding useful annual fraction $\fractionfmi{\stationindex}$. The map clearly indicates that wind power has a high availability in regions along the coastline and in the northern parts of Finland.

\begin{table}[htbp]
	\caption{Summary statistics for the useful fraction of $2021$ to run a dishwasher (see Figure \ref{fig:dishwasher}) only with wind power.}
	\begin{center}
		\begin{tabular}{|c|c|c|c|c|}
			\hline
			\textbf{Battery capacity} &  \textbf{Min fraction} &  \textbf{Max fraction}&  \textbf{Mean} &  \textbf{std} \\
			\hline
			200Wh& 0.14& 0.96& 0.66 & 0.19 \\
			\hline
			500Wh & 0.18& 0.97& 0.70 & 0.18  \\
			\hline
			800Wh  & 0.20& 0.97& 0.72 & 0.18  \\
			\hline
			1000Wh  & \textbf {0.21}& \textbf {0.97}& \textbf {0.72} & \textbf {0.17}  \\
			\hline
			1500Wh  &\textbf {0.21}&\textbf {0.97}& \textbf {0.72} & \textbf {0.17}  \\
			\hline
			2000Wh  & \textbf {0.21}& \textbf {0.97}& \textbf {0.72} & \textbf {0.17} \\
			\hline
		\end{tabular}
	\end{center}
	\label{tab:dishwasher}
\end{table}

\begin{figure}[htb]
   \center{\includegraphics[width=60mm]
	       {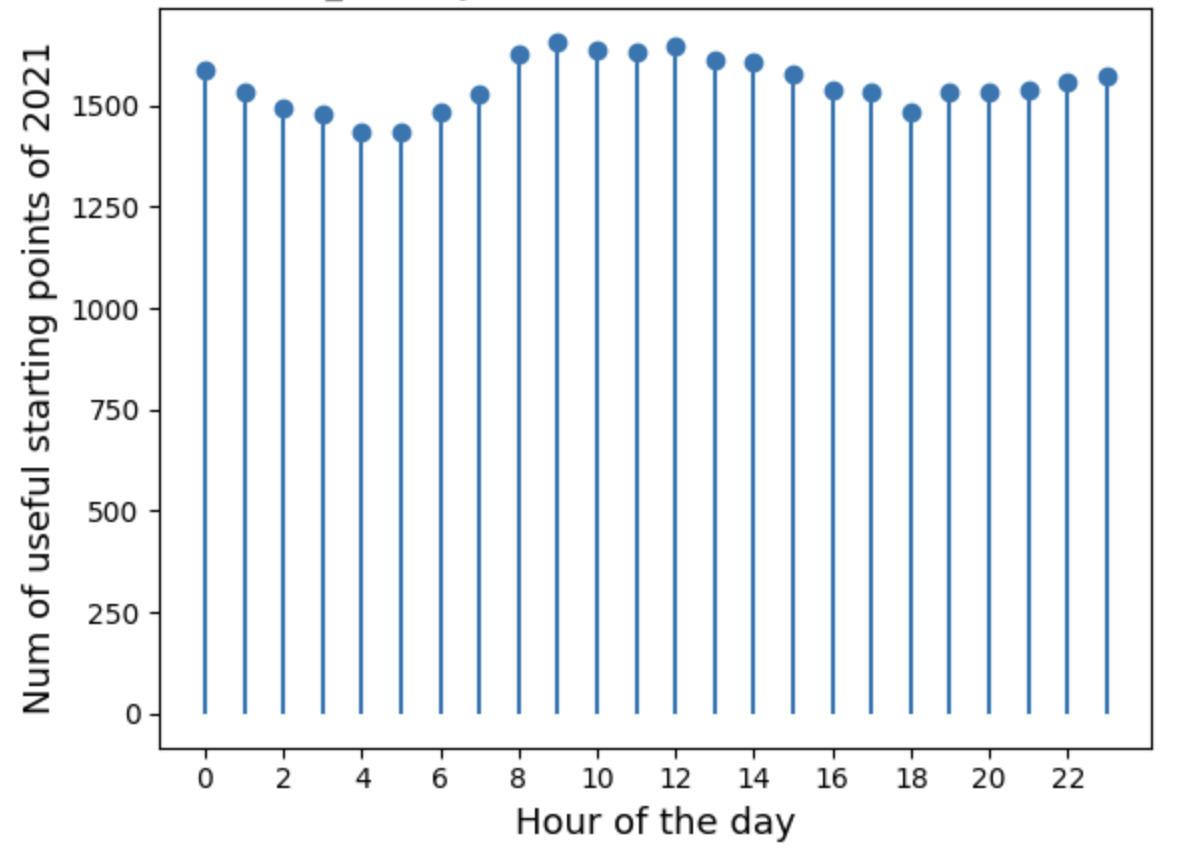}}
	 \caption{Distribution of useful starting points over 24h in Sotkamo Tuhkakyla weather station.}
	 \label{fig:distr}
\end{figure}

We also did further analysis to explore the periodic variation of useful fractions over 24 hours of the day and 12 months of the year. 
Figure \ref{fig:distr} shows an example weather station Sotkamo Tuhkakyla where useful annual 
fraction is 0.7, but the distribution over 24 hours is quite uniform. To quantize the uniformity of the distribution 
for each station, information entropy is applied; in particular, the following formula is used to calculate the 
entropy 
$$H=-\sum_{\houridx=1}^{24} p_\houridx \log p_\houridx$$ 
where $p_\houridx$ is the fraction of useful starting time points that fall into the $i$th hour of a day during 2021.

The result shows the entropies for all weather stations included in this study are larger than $0.97$, 
indicating the distribution of useful starting time points is quite uniform over the $24$ hours of a day. 

\begin{figure}
  \centering
     \begin{subfigure}[b]{0.23\textwidth}
         \centering
         \includegraphics[width=\textwidth]{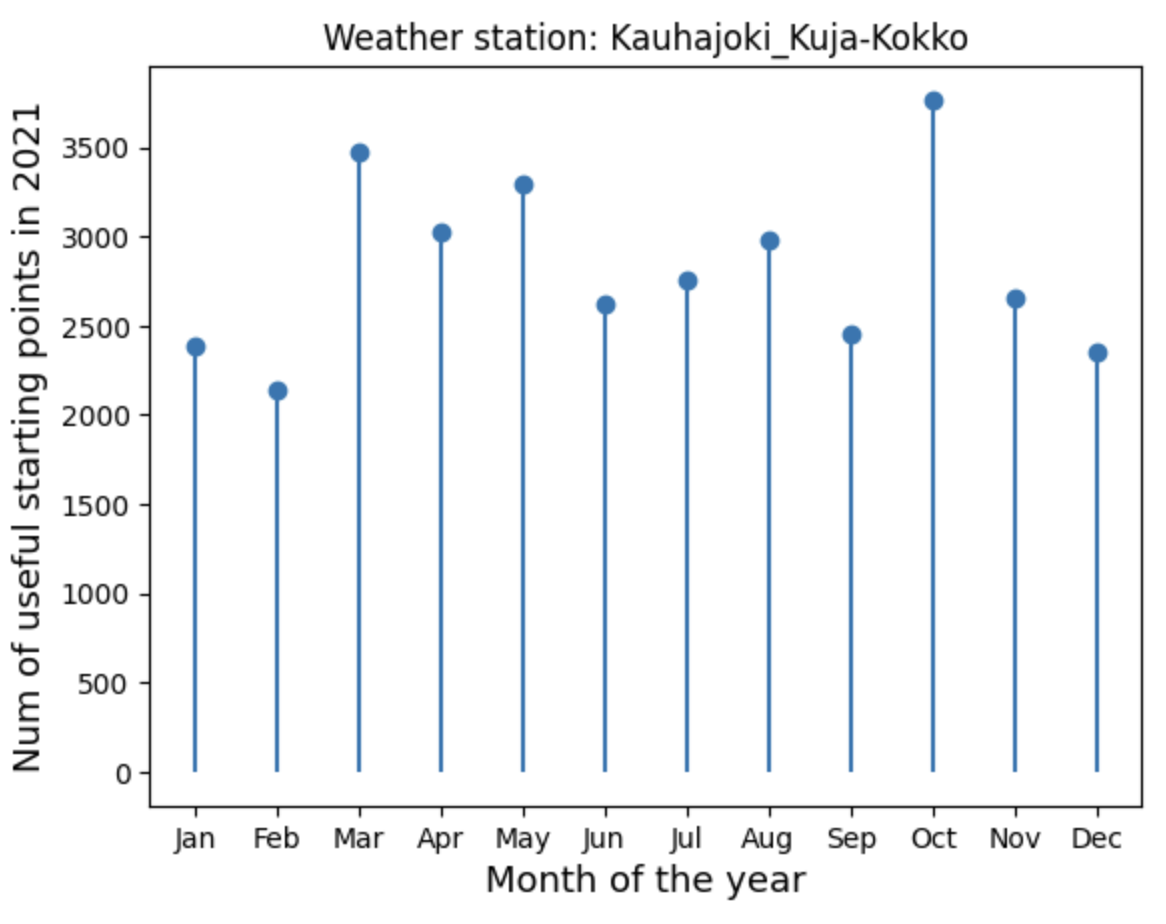}
         \caption{}
         \label{sub:d1}
     \end{subfigure}
       \begin{subfigure}[b]{0.24\textwidth}
         \centering
         \includegraphics[width=\textwidth]{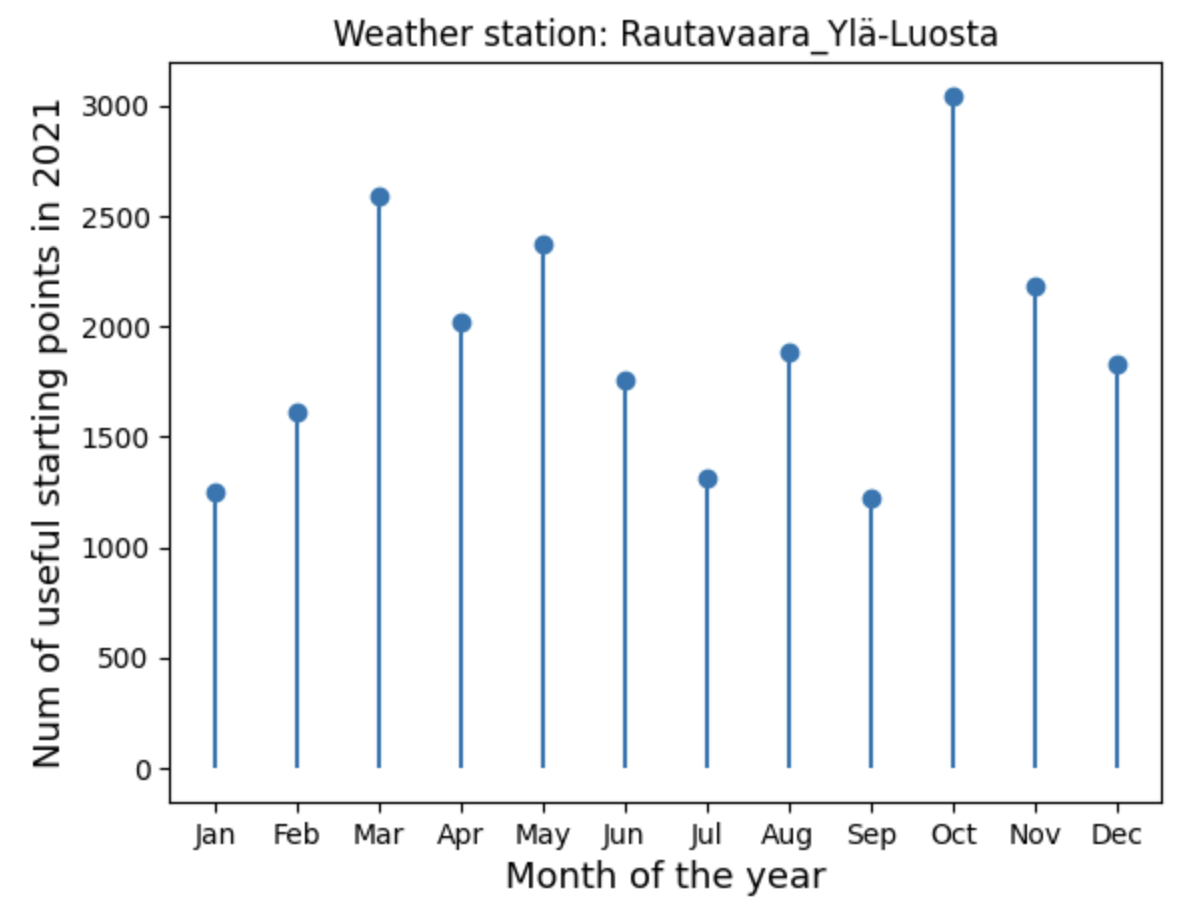}
         \caption{}
           \label{sub:d2}
     \end{subfigure}
       \begin{subfigure}[b]{0.23\textwidth}
         \centering
         \includegraphics[width=\textwidth]{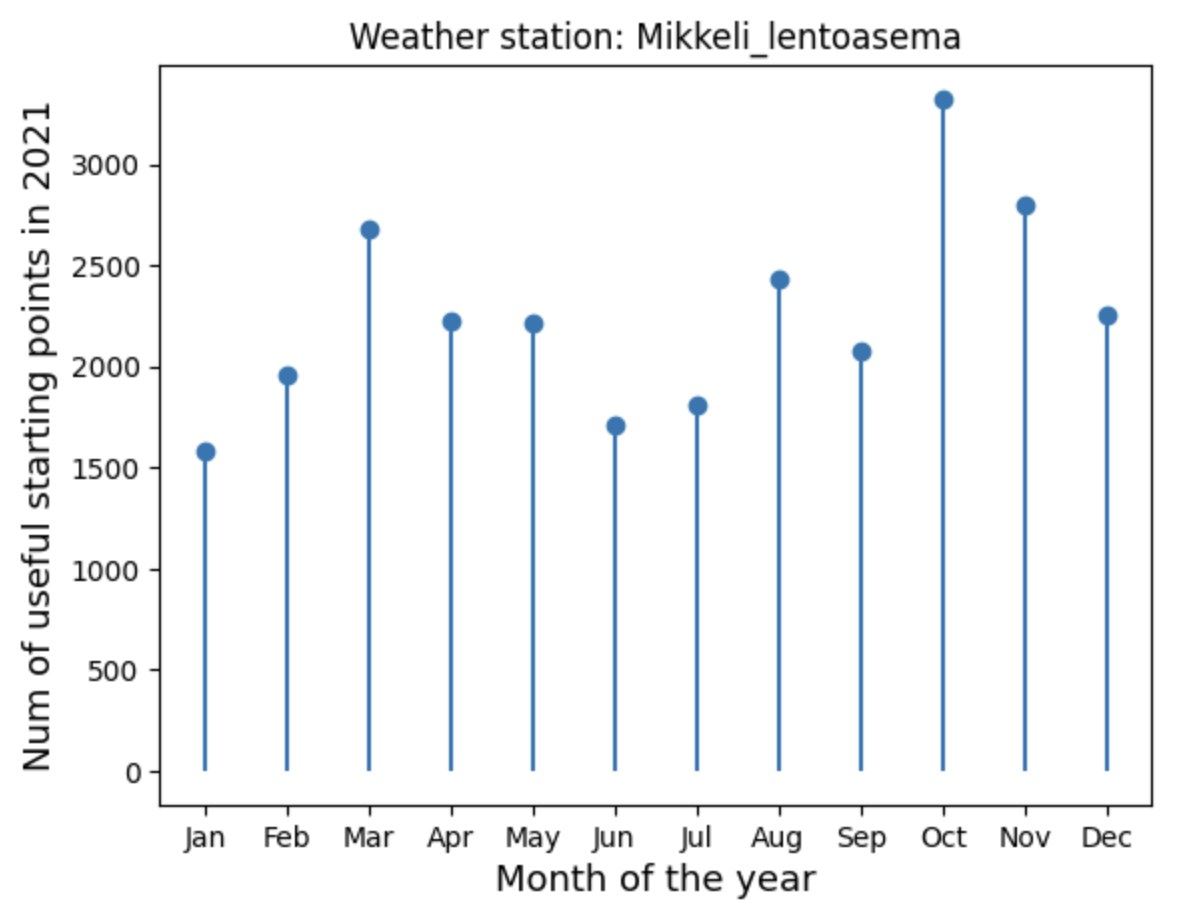}
         \caption{}
          \label{sub:d3}
     \end{subfigure}
       \begin{subfigure}[b]{0.24\textwidth}
         \centering
         \includegraphics[width=\textwidth]{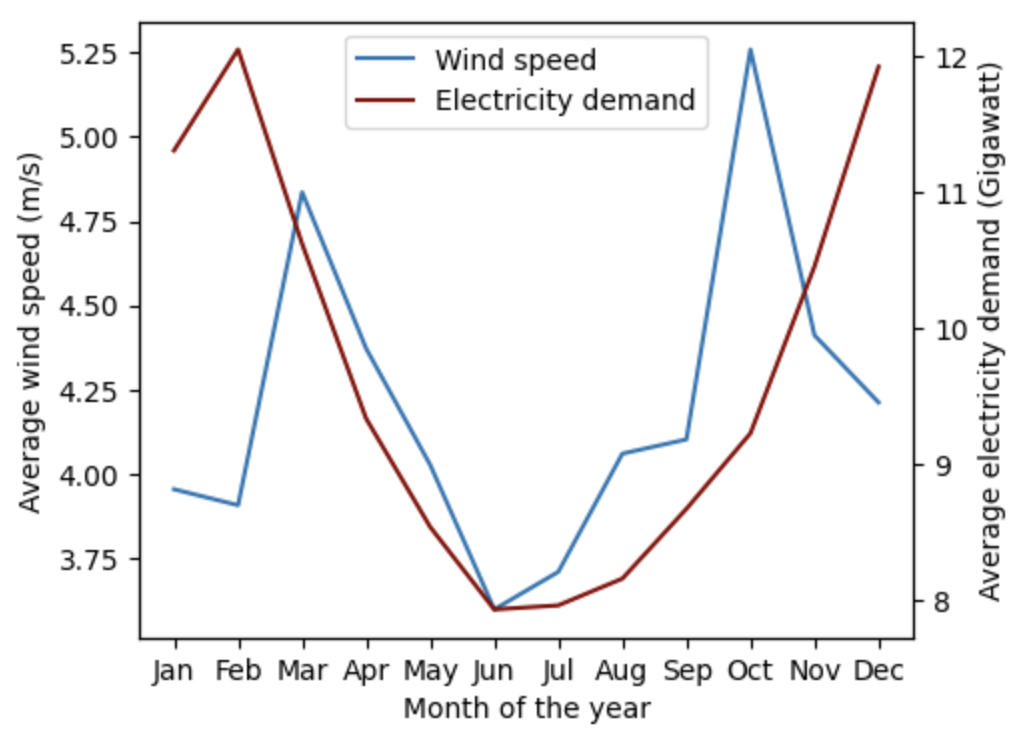}
         \caption{}
          \label{sub:dm}
     \end{subfigure}
    \caption{Distribution of useful starting points over 12 months at some weather stations (plots a,b,c) and a comparison between average wind speed and electricity demand of year 2021 (plot d)}
    \label{fig:months}
\end{figure}

Whereas the distribution characteristic over 12 months of the year is different. Figure \ref{fig:months} (\subref{sub:d1}, \subref{sub:d2}, \subref{sub:d3}) show significant seasonal trends in some weather stations,  and the trends match the general characteristics of wind speed in Finland: average wind speed of March and October is relatively higher (year 2021)(Figure \ref{fig:months} (\subref{sub:dm})). A further comparison with the trend of electricity demand in Finland (year 2021) shows during the summer months (June, July, and August), lower  electricity demand is aligned with  lower wind speed, but the trend is not quite consistent during winter months when the electricity demand achieves peaks. Electricity demand data is from FinGrid open database\cite{fingrid}.

%
%

\begin{figure}[htb]
   \center{\includegraphics[width=60mm]
	       {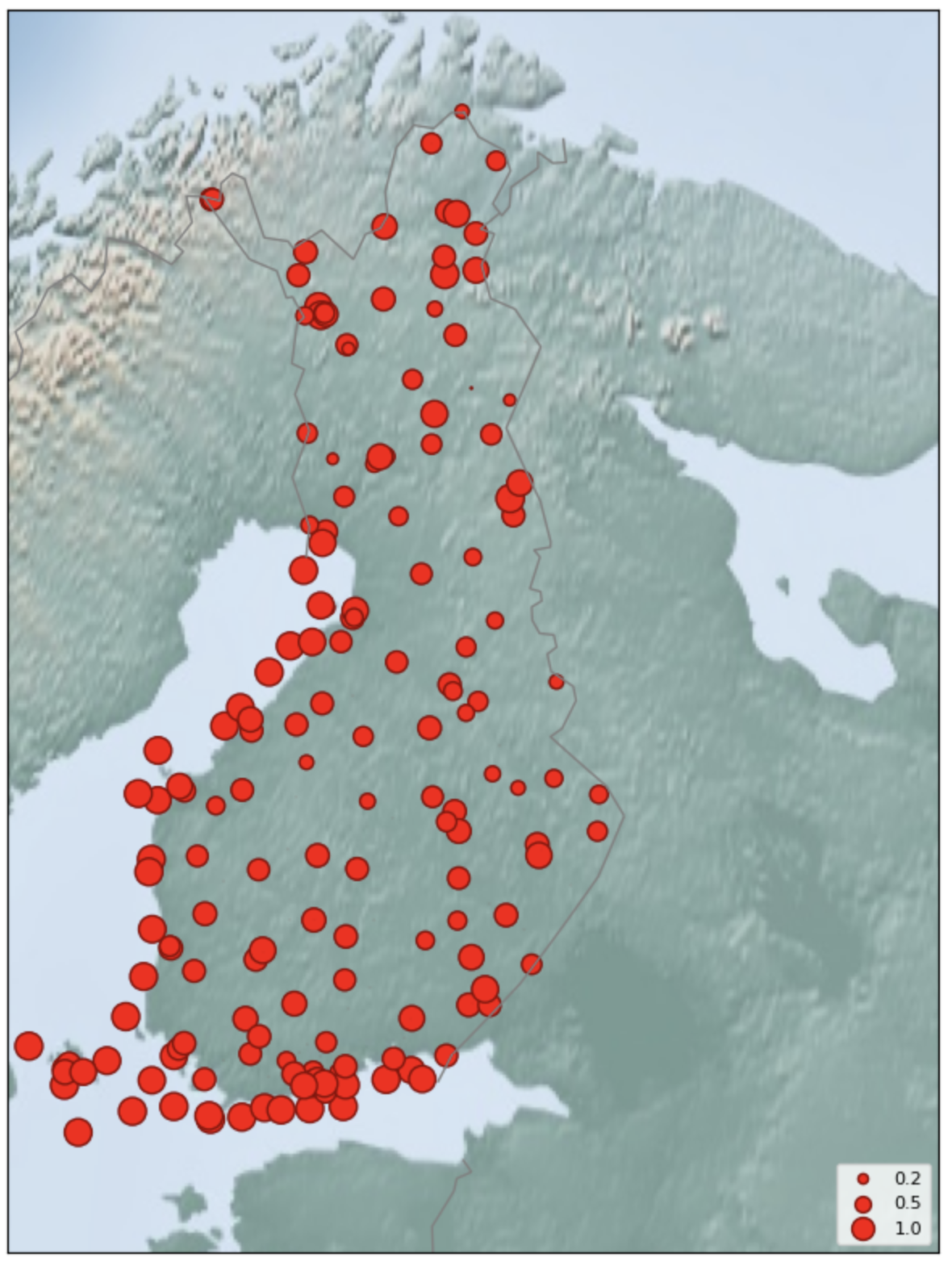}}

	 \caption{Spatial distribution of annual useful fractions $\fractionfmi{\stationindex}$ of $2021$ during which the 
	 	load profile $\household{\timeidx'}$ (see Figure \ref{fig:house}) could have been powered only from 
	 	wind power (using a battery capacity $\battcap=2500 {\rm Wh}$). The markers (red dots) represents FMI weather stations, 
	 	indexed by $\stationindex=1,\ldots,165$. The marker size (radius) is scaled by $\fractionfmi{\stationindex}$. }

	 \label{fig:heat_h}
\end{figure}

\begin{table}[htbp]
\caption{Useful annual fractions in FMI weather stations to provide wind power for an entire house, given different battery capacities. }
\begin{center}
\begin{tabular}{|c|c|c|c|c|}
   \hline
   \textbf{Battery capacity} &  \textbf{Min fraction} &  \textbf{Max fraction}&  \textbf{Mean} &  \textbf{std} \\
   \hline
      1000Wh& 0.17& 0.97& 0.68 & 0.18 \\
       \hline
      1500Wh & 0.19& 0.97& 0.70 & 0.18  \\
       \hline
       2000Wh  & 0.19& 0.97& 0.71 & 0.18  \\
       \hline
        2500Wh  & \textbf {0.20}& \textbf {0.97}& \textbf {0.72} & \textbf {0.17}  \\
         \hline
        3000Wh  &\textbf {0.20}&\textbf {0.97}& \textbf {0.72} & \textbf {0.17}  \\
         \hline
    \end{tabular}
  \end{center}
\label{tab:house}
\end{table}
From Table \ref{tab:house}, we can see compared to running a dishwasher, larger battery capacity is needed to fully exploit the wind resources to provide power for an entire house with profile $\household{\timeidx'}$. Similarly, a map (Figure \ref{fig:heat_h} ) is generated to visulise the useful annual fractions of 2021.

\begin{figure}[htb]
	\center{\includegraphics[width=80mm]
		{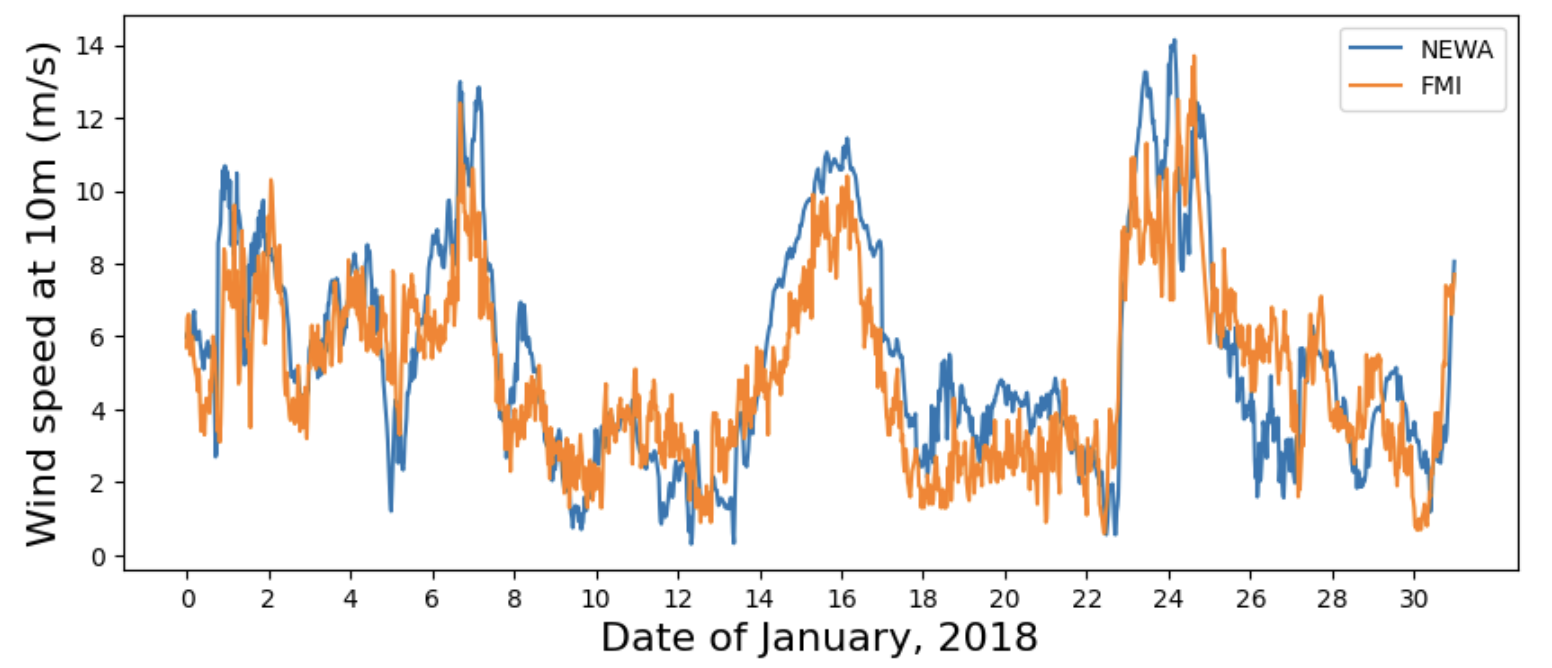}}
	\caption{Comparing modelled wind speed from NEWA with FMI wind speed observations during Year 2018 at Helsinki Kumpula) }
	\label{fig:newa}
\end{figure}

\begin{figure}[htb]
	\center{\includegraphics[width=80mm]
		{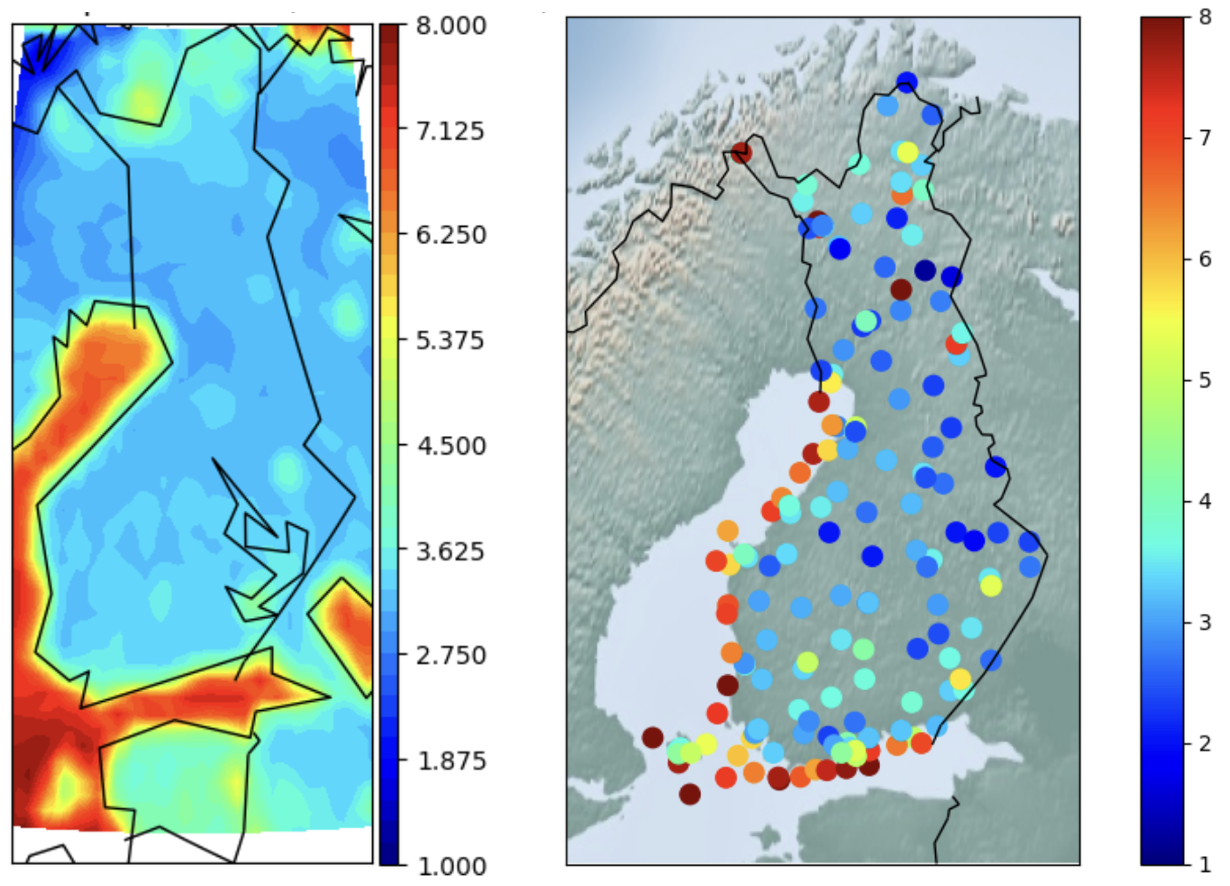}}
	
	\caption{Left: ERA5 reanalysis data for the average wind speed at $10$ m height during year $2021$. Right: Dots represent FMI 
		weather stations with dot color represents the average observed wind speed at $10$ m height during year $2021$.}
	\label{fig:era5}
\end{figure}

{\bf Comparison with NEWA and ERA5.} While wind speeds observations from FMI have a high temporal resolution (10-min intervals), they only allow for a poor spatial resolution of the resulting wind power atlas. Some sites are quite close to each other, with distances less than 1km, and some sites are far away from their neighbors, with distances larger than 50km, though generally, the sites are well distributed.
Thus, the wind atlas in Figure \ref{fig:heat} and \ref{fig:heat_h} are mostly useful for wind turbine locations nearby FMI 
weather stations. One possible extension of our approach would be to use high-resolution wind reanalysis data sets to 
interpolate between FMI weather stations. To this end, we compared FMI observation with reanalysis datasets
from NEWA\cite{newa} and ERA5\cite{era5}, which use the Weather Research and Forecasting Model (WRF) and the European Centre for Medium-Range Weather Forecasts Integrated Forecasting System (ECMWF-IFS), respectively. As NEWA only updated the reanalysis dataset to the year 2018, so we also downloaded corresponding data of the year 2018 from FMI for comparison. From Figure \ref{fig:newa} and \ref{fig:era5}, we can see NEWA data and ERA5 data generally compare well to FMI observations except for some northern FMI stations of Finland. 

\section{Conclusion}

We have used open weather data from the FMI to construct a wind power atlas for Finland. Unlike existing approaches for constructing a wind power atlas, our method employed a high temporal resolution that allows 
to verify whether a particular load profile could be powered solely by a wind turbine (combined with a battery storage capacity). Our results indicate that wind power has a high availability in regions along the 
coastline and in northern parts of Finland. We highlight that our analysis was exploratory using historical 
observation data from previous years. As an important next step, we will consider the short-term predictability 
of wind power in Finland. 

\bibliographystyle{IEEEtran}
\bibliography{IEEEabrv,Windpower.bib}

\end{document}